\documentclass[3p]{elsarticle}

\usepackage{amssymb}
\usepackage{amsmath}
\usepackage{theorem}

\newtheorem{theo}{Theorem}[section]
{\theorembodyfont{\rmfamily}
\newtheorem{defin}[theo]{Definition}
}
\newtheorem{prop}[theo]{Proposition}

{\theorembodyfont{\rmfamily}
\newtheorem{example}[theo]{Example}
}
{\theorembodyfont{\rmfamily}

}
{\theorembodyfont{\rmfamily}

}
\begin{document}

\title{Implementing generating functions to obtain power indices with coalition configuration}

\author[rv]{J. Rodr\'{\i}guez-Veiga}
\ead{jorge.rodriguez.veiga@gmail.com}

\author[nf]{G.I. Novoa-Flores}
\ead{gnovoaflores@gmail.com}

\author[cm]{B. Casas-M\'{e}ndez}
\ead{balbina.casas.mendez@usc.es}

\cortext[cor1]{Corresponding author: Jorge Rodr\'{\i}guez Veiga. ITMATI, R\'ua de Constantino Candeira s/n. Campus Vida. Santiago de Compostela (A Coru\~{n}a). Spain. Postal code: 15782. Phone: +34 650470360.}

\address[rv]{ITMATI (Technological Institute of Industrial Mathematics), Santiago de Compostela, Spain.}

\address[nf]{Department of Statistics and Operations Research, University of Santiago de Compostela, Spain.}

\address[cm]{MODESTYA Research Group, Department of Statistics and Operations Research, University of Santiago de Compostela, Spain.}

\begin{abstract}

We consider the Banzhaf-Coleman and Owen power indices for weighted majority games modified by a coalition configuration. We present calculation algorithms of them that make use of the method of generating functions. We programmed the procedure in the open language \texttt{R} and it is illustrated by a real life example taken from social sciences.

\end{abstract}


\maketitle

\bigskip

\section{Introduction}

In a cooperative game with transferable utility, also called TU-game, agents have mechanisms that allow them to make binding agreements. The utility derived from the cooperation of any set of agents (or players) can be transferred and divided in any way between them. One of the main objectives of cooperative game theory is the study of values, that is, distribution rules that assign to each game a vector, so that each coordinate vector represents the payment allocated to each of the players. There are numerous values, two of which, widely used, are the Shapley value (Shapley, \cite{sha53}) and the Banzhaf value (Banzhaf, \cite{ban65}).

An important family of TU-games is formed by simple games, which have interesting applications, especially in the field of political science. Within the family of simple games, weighted majority games play an important role. For example, a Parliament can be seen as a weighted majority game, assuming that players are the political parties and decisions are taken by majority.

When working with simple games, instead of talking about value, it is often used the term power index as simple games are usually used as models of decision-making bodies in which agreements are taken by voting. The interest in these games usually focuses on knowing the power or influence having a player in the final result. In this context, the Shapley value and the Banzhaf value are renamed Shapley-Shubik power index (Shapley and Shubik, \cite{sha54}) and Banzhaf-Coleman power index (Banzhaf, \cite{ban65} and Coleman, \cite{col71}), respectively. Numerous works in the literature are devoted to the definition, study, calculation, and applications of power indices.

The power indices can be obtained in an exact or approximated manner by using different tools. One of the most commonly used is given by the called generating functions. This, in the case of weighted majority games, is based on the use of a technique of combinatorial analysis. Roughly speaking, a generating function is a polynomial that allows to enumerate the set of possible coalitions, while having control about their respective weights. This is very useful because it allows to obtain the exact value of the indices, even in games with many players, making use of algorithms that can be programmed with a computer language. This technique was used by David G. Cantor (1962) (cf. Lucas, \cite{luc83}) for the Shapley-Shubik power index and by Brams and Affuso (\cite{bra76}) for calculating the Banzhaf-Coleman power index.

In a more general model of TU-games with coalition structure, a partition is considered over the set of players, because they could have some preferences to join others motivated by some reasons as ideological features, or geographic location. For this reason, the Owen value (Owen, \cite{owe77}) was proposed for TU-games with coalition structure, which generalizes the Shapley value, and the Banzhaf-Owen power index (Owen, \cite{owe82}) was proposed  for simple games with coalition structure as an extension of the Banzhaf-Coleman power index. Alonso-Meijide and Bowles (\cite{alo05}) used generating functions to calculate the Owen and Banzhaf-Owen values for weighted majority games.

Nonetheless, this model might not be appropiate in some situations in which players could prefer to join some players for some reasons and join some others for other reasons. Let us illustrate this situation with an example taken from Andjiga and Courtin (\cite{and13}): Consider the diplomatic relations among countries. In real life, countries are organized into international coalitions not necessarily disjoined. For example, France and Spain, among others, belong to the European Union and the North Atlantic Treaty Organization (NATO), while The United States belongs to NATO and NAFTA (a coalition with Mexico and Canada). For this reason, the more general model of games with coalition configuration was introduced. In this model, it is considered a cover of the set of players instead of a partition of it. Albizuri and Aurrekoetxea (\cite{alb06}) and Albizuri, Aurrekoetxea, and Zarzuelo (\cite{alb062}) propose the generalized Banzhaf-Coleman index and the configuration value\footnote{We use for these values the names proposed in Albizuri and Aurrekoetxea (\cite{alb06}) and Albizuri, Aurrekoetxea, and Zarzuelo (\cite{alb062}).}, for the models of simple games and TU-games, respectively, with a coalition configuration. These rules generalize the Banzhaf-Coleman power index and the Owen value, respectively.

In the current paper, we apply the method of generating functions to calculate the generalized Banzhaf-Coleman index and the configuration value for the family of weighted majority games with coalition configuration in such a way that we extend the results of Lucas (\cite{luc83}), Brams and Affuso (\cite{bra76}), and Alonso-Meijide and Bowles (\cite{alo05}). We illustrate the algorithms with a small numerical example and another one taken from real life. We also show an \texttt{R} code that implements the new algorithms.

We start the work with some preliminaries.

\section{Preliminaries}

A {\it cooperative game with transferable utility} (shortly a {\it TU-game}) is a pair
$(P, f)$, where $P = \{1, \dots, p\}$ is a {\it set of players} and $f$ is a real function (the {\it characteristic function}) which allocates to each
coalition $T \subseteq P$ a real value. This value can be interpreted as the payoff
obtained by each subset of players as the result of the cooperation among its members.

We say that a TU-game is {\it simple} if $f$ only takes the values 0 or 1, $f(P) = 1$, and the
game is monotonic, namely, $f(T) \leq f(S)$ if $T \subseteq S\subseteq P$. We represent for $W(f)$ the
set of all {\it winning coalitions} of the game $(P,f)$, it means, $W(f) = \{T \subseteq P\, :\, f(T) = 1\}$.
In this paper, we will work with a particular class of simple games, the {\it weigthed majority games}, in which  exists a set of {\it weights} $w_1,\dots,w_p$, with $w_i\in \Bbb{N}\cup \{0\}$ for $i = 1,\dots, p$, and a {\it quota} $q \in \mathbb{N}$ ($q>0$), such that for each $T\subseteq P$, $f(T) = 1$ if and only if $w(T) \geq q$ where $w(T) = \sum_{i \in T}w_{i}$. We will represent a weighted majority game by $\lbrack q; w_{1},  \dots, w_{p} \rbrack$.

Given a sequence of real numbers, $\{a_{n}\}_{n \in \mathbb{N}}$, we can associate the following series:
$$f_a(x)=\sum_{j \geq 0}a_jx^j.$$

This series is called the {\it generating function} of the sequence $a_{n}$, and may be finite or infinite.

In this paper, we will work with generating functions of several variables, which can be written as $f_a(x,y,z)=\sum_{j \geq 0}\sum_{k \geq 0}\sum_{l \geq 0}a_{j,k,l}x^jy^kz^l$
where $a_{j,k,l}$ represents a number depending on $j$, $k$ and $l$.

Let us consider a set of players $P = \{1, \dots, p\}$. A {\it coalition configuration}  $C = \{C_{1}, \dots, C_{c}\}$
is a cover of $P$, it is, a collection of non-empty subsets of $P$ such that $\bigcup_{k=1}^{c}C_{k}=P$. As a consequence of the fact that each player
could be simultaneously in more than one element of $C$, for each $i\in P,$ we will consider $C^{i} =
\{C_{k} \in C\, :\, i \in C_{k}\}$ as the set of elements of $C$ in which player $i$ is in. A {\it TU-game with coalition configuration} is a triple $(P,f;C)$ where $(P,f)$ is a TU-game and $C$ a coalition configuration over $P$. Let us denote by $SC(P)$ the family of simple games with coalition configuration and players set $P$.

Albizuri and Aurrekoetxea (\cite{alb06}) generalize the Banzhaf-Coleman index to simple games with coalition configuration, which we call here {\it generalized Banzhaf-Coleman index}, and it is defined for each simple game with coalition configuration $(P,f;C)$ and each player $i \in P$ as:

\begin{equation}
\beta_{i}(P,f;C)=\sum\limits_{C_{k} \in C^{i}}\sum\limits_{R \subseteq C\backslash C^{i}}\sum\limits_{T \subseteq C_{k}\backslash \{i\}}\frac{1}{2^{c+c_{k}-2}}\Bigl(f\bigl(T \cup S_R \cup \{i\}\bigr)- f(T \cup S_R)\Bigr),
\end{equation}

\noindent with $c_k=\mid C_k \mid$ and $S_R=\cup_{C_l\in R}C_l.$

On the other hand, Albizuri, Aurrekoetxea, and Zarzuelo (\cite{alb062}) generalize the Owen value to games with coalition configuration, which we call here {\it configuration index} when it is applied to simple games with coalition configuration and it is defined for each simple game with coalition configuration $(P,f;C)$ and each player $i \in P$ as:

\begin{equation}
\Phi_{i}(P,f;C)=\displaystyle{\sum\limits_{C_{k} \in C^{i}}\sum\limits_{R \subseteq C\backslash C^{i}}\sum\limits_{T \subseteq C_{k}\backslash \{i\}}q(c_k,r,t)}\Bigl(f\bigl(T \cup S_R \cup \{i\}\bigr)- f(T \cup S_R)\Bigr),
\end{equation}

\noindent with $q(c_k,r,t)=\dfrac{r!(c-r-1)!}{c!}\dfrac{t!(c_k-t-1)!}{c_k!}$, $r=\mid R\mid$, and $t=\mid T\mid$.

\section{Generating functions for indices with coalition configuration}

In this section, we present the procedures based on the generating functions for calculating the two power indices introduced in the previous section, the generalized Banzhaf-Coleman index and the configuration index.

\begin{defin}

Let us consider $(P,f;C)\in SC(P)$ and $i\in P$. We say that $T \subseteq P$ is a {\it coalition consistent with the coalition configuration for $C_{k} \in C^{i}$}, if $T=\left(\cup_{C_{l} \in R}C_{l}\right)\cup S$ with
$R\subseteq C\backslash C^{i}$, and $S \subseteq C_{k}$.

\end{defin}

The set of all coalitions consistent with the coalition configuration for $C_{k} \in C^{i}$ will be denoted by ${\cal C}(i,C_k,C)$. According to this definition, both indices, the generalized Banzhaf-Coleman index and the configuration index, for each player, are weighted sums of the contribution of such player to the various coalitions consistent with coalitions (to which he belongs to) of the configuration to which he can join. In the first index, the weights do not consider the size of those coalitions consistent while the second index does. As it is known, we can give interpretations with the same flavor to the rules of Banzhaf, Shapley, Banzhaf-Owen, and Owen.

\begin{defin}

Let us consider $(P,f;C)\in SC(P)$ and $i\in P$. We say that $T\subseteq P$ is a {\it significant coalition consistent with the coalition configuration for $C_{k} \in C^{i}$}, if $T \notin W(f)$, $T \cup \{i\} \in W(f)$, and $T \in {\cal C}(i,C_k, C)$.

\end{defin}

The {\it number of all significant coalitions consistent with the coalition configuration for $C_{k} \in C^{i}$} will be denoted by $\sigma^{i,k}(P,f;C)$.

\subsection{The generalized Banzhaf-Coleman index}

For each simple game with coalition configuration $(P,f;C)$ and each player $i \in P$, from (1) we deduce that:

$$
\beta_{i}(P,f;C)=\sum\limits_{C_{k} \in C^{i}}\frac{\sigma^{i,k}(P,f;C)}{2^{c+c_{k}-2}}.
$$

The following result allows us to efficiently compute the number of significant coalitions consistent with the coalition configuration for $C_{k} \in C^{i}$ in a weighted majority game.

\begin{prop}

Let $(P,f;C)$ be a weighted majority game with coalition configuration, where the game is represented by $[q; w_{1}, \dots, w_{p}]$ and $C=\{C_{1}, \dots, C_c\}$. It holds that:

$$
\sigma^{i,k}(P,f;C)=\sum\limits_{m=q-w_i}^{q-1}\nu^{i,k}_m,
$$	

\noindent being $\nu^{i,k}_m=|\{T \in {\cal C}(i,C_k,C)\, :\,  i \notin T \hspace{0.15cm} {\rm and} \hspace{0.15cm} w(T)=m\}|$ for each player $i \in P$, each coalition $C_k \in C^i$, and $m\in \Bbb{N}$.
\end{prop}

\noindent {\bf Proof.} Let us consider $(P,f;C)$ a weighted majority game with coalition configuration, where the game is represented by $[q; w_{1}, \dots, w_{p}]$, $C=\{C_{1}, \dots, C_c\}$, and take a player $i \in P$ and $C_k \in C^i$.

The coalitions consistent with the coalition configuration for $C_{k} \in C^{i}$ in which player $i$ is not in and with a weight between $q-w_i$ and $q-1$ are losing coalitions.	If player $i$ joins to any of those coalitions, the weight of the coalition is, at least, $q$ so it turns into winning coalitions, that is, they are significant coalitions consistent with the coalition configuration for $C_{k} \in C^{i}$. If we add the number of all such coalitions we obtain $\sigma^{i,k}(P,f;C)$. $\Box$\\

Now, we present the main result of this subsection in which a generating function to compute the numbers $\nu_{m}^{i,k}$ is proposed.

\begin{theo}\label{ThCC}
Let $(P,f;C)$ be a weighted majority game with coalition configuration, where the game is represented by $[q; w_{1}, \dots, w_{p}]$ and $C=\{C_{1}, \dots, C_c\}$. Given $i \in P$ and $C_k \in C^i$, the generating function of the numbers $\{\nu^{i,k}_m\}_{m \geq 0}$ is obtained from the function:

$$
G_{i,k}(x,(t_{j})_{j \in C_k \backslash \{i\} \cup \{C_l\, :\, C_l \notin C^i  \}})=\prod\limits_{{j} \in C_k \backslash \{i\}}(1+x^{w_{j}}t_{j})\prod\limits_{C_l \notin C^i}(1+x^{w(C_l)}[\prod_{j \in C_l}t_{j}]).
$$

Once the above function is computed, in each monomial and for each term of the form $t^{e_{j}}_{j}$,
we multiply by $x^{-(e_{j}-1)w_{j}}$, and then we remove all the terms of the form $t^{e_{j}}_{j}$. The resulting function, we call it $G_{i,k}(x)$, is the generating function of the numbers $\{\nu^{i,k}_m\}_{m \geq 0}$.
\end{theo}

\noindent {\bf Proof.} Let $(P,f;C)$ be a weighted majority game represented by $[q;w_1, \dots, w_p]$ with coalition configuration $C=\{C_1, \dots, C_c \}$, and consider a player $i \in C_k$ for some $C_k \in C^i$.

Let us consider the function:

\begin{equation}
\prod\limits_{{j} \in C_{k}}(1+x^{w_{j}}t_{j})
\prod\limits_{C_l \notin C^i}(1+x^{w(C_l)}[\prod_{j \in C_l}t_{j}]),
\end{equation}
\noindent and let us denote by $\hat{G}_{i,k}(x)$ the function obtained of the above one after multiplying each monomial by $x^{-(e_{j}-1)w_{j}}$, for each coefficient of the form $t_j^{e_j}$, and after removal such coefficients. It is clear that:

\bigskip

\begin{center}

\renewcommand*{\arraystretch}{1.75}{
\begin{tabular}{rl}
	$\hat{G}_{i,k}(x)$ =&
	$\sum\limits_{T \subseteq C_k}\sum\limits_{R \subseteq C\backslash C^i} x^{w(T \cup R)}
	=\sum\limits_{T \subseteq C_k}\sum\limits_{R \subseteq C\backslash C^i}x^{\sum\limits_{j \in T \cup R}w_j}
	=\sum\limits_{T \subseteq C_k}\sum\limits_{R \subseteq C\backslash C^i}\prod\limits_{j \in T\cup R}x^{w_j}$	    \\
	=&
	$
	1+\sum\limits_{S \in {\cal C}(i,C_k,C)}\prod\limits_{j \in S}x^{w_j}=1+\sum\limits_{S \in {\cal C}(i,C_k,C)}x^{\sum\limits_{j \in S}w_j}
	=1+\sum\limits_{S \in {\cal C}(i,C_k,C)}x^{w(S)}
	=\sum\limits_{m=0}^{w(P)}\nu_m^kx^m$,\\
	
\end{tabular}}
\bigskip

\end{center}

\noindent where $\nu_0^k=1$ and for $m>0$, $\nu_{m}^{k}=|\{S \in {\cal C}(i,C_k,C)\, :\, w(S)=m\}|$, where the last expression is obtained doing $w(S)=m$ and by grouping exponents of the same value.

Finally, in order to obtain the numbers $\{\nu^{i,k}_m\}_{m \geq 0}$, we delete the factor $(1+t_{j}x^{w_{j}})$ with $j=i$ in the equation $(3)$.
$\Box$\\

Theorem 3.4 provides a method for obtaining the numbers $\{\nu^{i,k}_m\}_{m\geq 0}$. But for the generalized Banzhaf-Coleman index, we need the numbers $\sigma^{i,k}(P,f;C)$. These values can be obtained by selecting in $G_{i,k}(x)$ the coefficients of the monomials $x^m$ in where the exponent $m$ of $x$ takes values between $q-w_i$ and $q-1$.

The following example will illustrate how to proceed to compute the generalized Banzhaf-Coleman index.

\begin{example}
Let us consider $P=\{1,2,3,4\}$ as the set of players, $C=\{C_1=\{1,2,3\},C_2=\{2,3\},C_3=\{3,4\}\}$ as the coalition configuration on $P$, and a weighted majority game represented by $[3;1,2,2,1]$.

Take player $1$ and the coalition $C_{1}$. We will build the function $G_{1,1}(x)$. Following Theorem 3.4:
$$
\begin{array}{rl}
	G_{1,1}(x,t_2,t_3,t_4) = & \prod_{{j} \in C_1 \backslash \{1\}}(1+x^{w_{j}}t_{j})\prod_{C_l \notin C^1}(1+x^{w(C_l)}[\prod_{j \in C_l}t_{j}]) \\
&\\
	= & \Bigl[(1+x^{w_2}t_2)(1+x^{w_3}t_3)\Bigr] \Bigl[\bigl(1+x^{w(C_2)}t_2t_3\bigr)\bigl(1+x^{w(C_3)}t_3t_4\bigr)\Bigr]\\
&\\
	= & 1+x^{w_2}t_2+x^{w_3}t_3+x^{w_2+w_3}t_2t_3 +x^{w(C_2)}t_2t_3+x^{w_2+w(C_2)}t^2_2t_3+x^{w_3+w(C_2)}t_2t^2_3\\
&\\
    + & x^{w_2+w_3+w(C_2)}t^2_2t^2_3+x^{w(C_3)}t_3t_4+x^{w_2+w(C_3)}t_2t_3t_4+x^{w_3+w(C_3)}t^2_3t_4\\
    &\\
	+ & x^{w_2+w_3+w(C_3)}t_2t^2_3t_4+x^{w(C_2)+w(C_3)}t_2t_3^2t_4+x^{w_2+w(C_2)+w(C_3)}t_2^2t_3^2t_4\\
&\\
	+ & x^{w_3+w(C_2)+w(C_3)}t_2t_3^3t_4+x^{w_2+w_3+w(C_2)+w(C_3)}t_2^2t_3^3t_4.\\
\end{array}
$$
\bigskip

In the last sum, each term represents a coalition consistent with the coalition configuration for $C_1$ that not contains to the player 1. The term of the form $t_j$ indicates the players involved and the exponents of $t_j$ indicate the number of times that the corresponding players appear in each subcoalition of the consistent coalition (a subcoalition is a subset of $C_1$ or one coalition of $C\backslash C_1$). The exponent of $x$ is the sum of the weights of the subcoalitions.

Now, we multiply each monomial by terms of the form $x^{-(e_{j}-1)w_{j}}$ being $e_{j}$ the exponent of each $t_{j}$ on the coefficient of the monomial. The resulting function is: \bigskip

\begin{center}

\renewcommand*{\arraystretch}{1.75}{
\begin{tabular}{rl}
	
	1+&$x^{w_2}t_2+x^{w_3}t_3+x^{w_2+w_3}t_2t_3 +x^{w(C_2)}t_2t_3+x^{w(C_2)}t^2_2t_3$\\
	+&$x^{w(C_2)}t_2t^2_3+x^{w(C_2)}t^2_2t^2_3+x^{w(C_3)}t_3t_4+x^{w_2+w(C_3)}t_2t_3t_4+x^{w(C_3)}t^2_3t_4$\\
	+&$x^{w_2+w(C_3)}t_2t^2_3t_4+x^{w_2+w(C_3)}t_2t_3^2t_4+x^{w_2+w(C_3)}t_2^2t_3^2t_4+x^{w_2+w(C_3)}t_2t_3^3t_4
	+x^{w_2+w(C_3)}t_2^2t_3^3t_4$.\\
\end{tabular}}
\end{center}
\bigskip

Then, in each term, the exponent of $x$ is the weight of a coalition consistent.
If we delete all the coeficients of the form $t_j$ in the above function we obtain:

\bigskip
\begin{center}
\begin{tabular}{rl}

1+& $x^{w_2}+x^{w_3}+x^{w_2+w_3}+x^{w(C_2)}+x^{w(C_2)}+x^{w(C_2)}+x^{w(C_2)}+x^{w(C_3)}+x^{w_2+w(C_3)}$\\
+&$x^{w(C_3)}+x^{w_2+w(C_3)}+x^{w_2+w(C_3)}+x^{w_2+w(C_3)}+x^{w_2+w(C_3)}+x^{w_2+w(C_3)}$.\\
\end{tabular}
\end{center}
\bigskip

This is precisely the generating function of the numbers $\{\nu^{1,1}_m\}_{m \geq 0}$, when we replace the weights of each player and the weights of each coalition:

\bigskip
\begin{center}
\begin{tabular}{rl}

$G_{1,1}(x)=$&$1+x^2+x^2+x^4+x^4+x^4+x^4+x^4+x^3+x^5+x^3+x^5+x^5$\\
 + & $x^5+x^5+x^5=1+2x^2+2x^3+5x^4+6x^5$.\\

\end{tabular}
\end{center}
\bigskip
Finally, we choose the coefficients of the terms whose exponent of $x$ is between $q-w_1=2$ and $q-1=2$, and we have that $\sigma^{1,1}(P,f;C)=2$. Since player 1 is only in one element of the coalition configuration, we have that:
$$
\beta_{1}(P,f;C)=\sum\limits_{C_k \in C^1}\frac{\sigma^{1,k}(P,f;C)}{2^{c+c_k-2}}=\frac{\sigma^{1,1}(P,f;C)}{2^{c+c_1-2}}=\frac{2}{2^{3+3-2}}=0.125.
$$

Proceeding analogously with the remaining players, we obtain that the generalized Banzhaf-Coleman index for this game is $\beta(P,f;C)=\bigl(0.125,0.25,0.375,0.125\bigr).$ The same result is obtained by using the code provided in the Appendix. In this case, we introduce the function \texttt{BanzhafWeightedMajorityCC} and then we write:

\bigskip
\begin{center}
\begin{tabular}{l}

$\texttt{wi<-c(1,2,2,1)}$\\
$\texttt{q<-3}$\\
$\texttt{C<-list(c(1,2,3),c(2,3),c(3,4))}$\\
$\texttt{BanzhafWeightedMajorityCC(wi,q,C)}\  \Box$\\

\end{tabular}
\end{center}
\bigskip

\end{example}

\subsection{The configuration index}

For each simple game with coalition configuration $(P,f;C)$ and each player $i \in P$, from (2) we deduce that:\\
\begin{center}
\resizebox*{\linewidth}{5cm}{
	\begin{minipage}{\linewidth}
		\begin{align*}
		\Phi_{i}(P,f;C) &= \displaystyle {\sum\limits_{C_{k} \in C^{i}}\sum\limits_{\stackrel{T \in {\cal C}(i,C_k,C)}{\stackrel{T\not\in W(f)}{T\cup \{i\}\in W(f)}}}
			\frac{\mid {\overline C}^i(T)\mid!(\mid C\mid-\mid {\overline C}^i(T)\mid-1)!}{\mid C\mid!}}
		\displaystyle{\frac{\mid T \cap C_k \mid!(\mid C_k \mid- \mid T \cap C_k \mid -1)!}{\mid C_k\mid!}} \\
		& \\
		&=	\displaystyle{\sum_{C_{k} \in C^{i}}}
		\displaystyle{\sum_{r=0}^{c-c^i }}
		\displaystyle{\sum_{t=0}^{c_k -1}}
		\frac{r!(c-r-1)!}{c!}\frac{t!(c_k-t-1)!}{c_k!}\sigma^{i,k}_{r,t}(P,f;C),
		 \\
		\end{align*}
	\end{minipage}
}
\bigskip
\end{center}
\noindent where ${\overline C}^i(T)= \{ C_l \in C \backslash C^i\, : \, C_l\subseteq T \}$, $c^i=|C^i|$, and $\sigma^{i,k}_{r,t}(P,f;C)$ is the number of all significant coalitions consistent with the coalition configuration for $C_k\in C_i$, which contain $r$ coalitions of $C\backslash C^i$ and $t$ players belonging to $C_k$.

Let us see below how we  compute efficiently the total number of such coalitions. The proof is left to the reader. As in Proposition 3.3, let $(P,f;C)$ be a weighted majority game with coalition configuration, where the game is represented by $[q; w_{1}, \dots, w_{p}]$ and $C=\{C_{1}, \dots, C_c\}$. It holds that:

$$
\sigma^{i,k}_{r,t}(P,f;C)=\sum\limits_{m=q-w_i}^{q-1}\nu^{i,k}_{m,r,t},
$$	
\noindent being $\nu^{i,k}_{m,r,t}=|\{S \in {\cal C}(i,C_k,C)\, :\,  i \notin S, \hspace{0.15cm}w(S)=m, \hspace{0.15cm} \mid {\overline C}^i(S)\mid=r, \hspace{0.15cm} {\rm and} \hspace{0.15cm} \mid S\cap C_k\mid =t\}|,$ for each player $i \in P$, each coalition $C_k \in C^i$, $m\in \Bbb{N}$, $r\in \Bbb{N}$, and $t\in \Bbb{N}$.

Now, we present the main result of this subsection in which a generating function to compute the numbers $\nu_{m,r,t}^{i,k}$ is proposed.

\begin{theo}\label{ThCC2}
Let $(P,f;C)$ be a weighted majority game with coalition configuration, where the game is represented by $[q; w_{1}, \dots, w_{p}]$ and $C=\{C_{1}, \dots, C_c\}$. Given $i \in P$ and $C_k \in C^i$, the generating function of the numbers $\{\nu^{i,k}_{m,r,t}\}_{m \geq 0,\, r\geq 0,\,t\geq 0}$ is obtained from the function:

$$G'_{i,k}(x,(t_{j})_{j \in C_k \backslash \{i\} \cup \{C_l\, :\, C_l \notin C^i  \}},u,v)=
\prod\limits_{{j} \in C_k \backslash \{i\}}(1+x^{w_{j}}t_{j}u)
\prod\limits_{C_l \notin C^i}(1+x^{w(C_l)}[\prod_{j \in C_l}t_{j}]v).
$$

Once the above function is computed, in each monomial and for each term of the form $t^{e_{j}}_{j}$,
we multiply by $x^{-(e_{j}-1)w_{j}}$, and then we remove all the terms of the form $t^{e_{j}}_{j}$. The resulting function, we call it $G'_{i,k}(x,u,v)$, is the generating function of the numbers $\{\nu^{i,k}_{m,r,t}\}_{m \geq 0,\, r\geq 0,\,t\geq 0}$.
\end{theo}

\noindent {\bf Proof.} Let $(P,f;C)$ be a weighted majority game represented by $[q;w_1, \dots, w_p]$ with coalition configuration $C=\{C_1, \dots, C_c \}$, and consider a player $i \in C_k$ for some $C_k \in C^i$.

Let us consider the function:

\begin{equation}
\prod\limits_{{j} \in C_{k}}(1+x^{w_{j}}t_{j}u)
\prod\limits_{C_l \notin C^i}(1+x^{w(C_l)}[\prod_{j \in C_l}t_{j}]v),
\end{equation}

\noindent and let us denote by $\hat{G}'_{i,k}(x,u,v)$ the function obtained of the above one after multiplying each monomial by $x^{-(e_{j}-1)w_{j}}$, for each coefficient of the form $t_j^{e_j}$, and after removal such coefficients. It is clear that:

\bigskip
\begin{center}
\begin{tabular}{rl}
$\hat{G}'_{i,k}(x,u,v)$ =&
$\sum\limits_{T \subseteq C_k}\sum\limits_{R \subseteq C\backslash C^i} x^{w(T \cup R)}u^{t} v^{r}$\\
=& $\sum\limits_{T \subseteq C_k}\sum\limits_{R \subseteq C\backslash C^i}
x^{\sum\limits_{j \in T \cup R}w_j}u^{t} v^{r}$	    \\
=&
$\sum\limits_{T \subseteq C_k}\sum\limits_{R \subseteq C\backslash C^i}
\prod\limits_{j \in T\cup R}x^{w_j}u^{t} v^{r}$\\
=& $1+\sum\limits_{S \in {\cal C}(i,C_k,C)} \prod\limits_{j \in S\cap C_k}[x^{w_j}u]\prod\limits_{C_l \in {\overline C}^i(S)}[x^{w(C_l)}v] $\\
=&$1+\sum\limits_{S \in {\cal C}(i,C_k,C)}x^{\sum\limits_{j \in S}w_j}u^{\mid S\cap C_k\mid}v^{\mid {\overline C}^i(S)\mid}$\\=
& $1+\sum\limits_{S \in {\cal C}(i,C_k,C)}x^{w(S)}u^{\mid S\cap C_k\mid}v^{\mid {\overline C}^i(S)\mid} $\\
=& $\sum\limits_{m=0}^{w(P)}\sum\limits_{r=0}^{c-c^i}\sum\limits_{t=0}^{c_k}\nu_{m,r,t}^kx^mu^tv^r$,\\

\end{tabular}
\end{center}
\bigskip

\noindent where $\nu_{0,0,0}^k=1$ and for $m>0$, $\nu_{m,r,t}^{k}=|\{S \in {\cal C}(i,C_k,C)\, :\, w(S)=m,\, \mid {\overline C}^i(S)\mid = r,\, {\rm and}\, \mid S\cap C_k\mid= t\}|$, where the last expression is obtained doing $w(S)=m$, $\mid {\overline C}^i(S)\mid = r$, $\mid S\cap C_k\mid= t$, and by grouping exponents of the same value.

Finally, in order to obtain the numbers $\{\nu^{i,k}_{m,r,t}\}_{m \geq 0,\, r\geq 0,\, t\geq 0}$, we delete the factor $(1+x^{w_{j}}t_{j}u)$ with $j=i$ in the equation $(4)$.
$\Box$\\

Theorem 3.6 provides a method for obtaining the numbers $\{\nu^{i,k}_{m,r,t}\}_{m\geq 0,r\geq 0,t\geq 0}$. But for the configuration index, we need the numbers $\sigma^{i,k}_{r,t}(P,f;C).$ These values $\sigma^{i,k}_{r,t}(P,f;C)$ can be identified with a polynomial of the form:

$$g_{i,k}(u,v)=\sum_{r=0}^{c-c^i}\sum_{t=0}^{c_k-1}\sigma^{i,k}_{r,t}(P,f;C)u^tv^r,$$

\noindent and as $\sigma^{i,k}_{r,t}(P,f;C)=\sum_{m=q-w_i}^{q-1}\nu^{i,k}_{m,r,t}$, we have that:

$$g_{i,k}(u,v)=\sum_{r=0}^{c-c^i}\sum_{t=0}^{c_k-1}
\Bigl[\sum_{m=q-w_i}^{q-1}\nu^{i,k}_{m,r,t}\Bigr]u^tv^r.$$

Now, by Theorem 3.6,

$$G'_{i,k}(x,u,v)=
\sum_{r=0}^{c-c^i}\sum_{t=0}^{c_k-1}
\Bigl[\sum_{m=0}^{w(P)-w_i}\nu^{i,k}_{m,r,t}x^m\Bigr]u^tv^r,$$

\noindent from which it follows that to determine the coefficients of $g_{i,k}(u,v)$, we have to select the coefficients of the monomials $x^mu^tv^r$ in where the exponent $m$ of $x$ takes values between $q-w_i$ and $q-1$,  for each pair of exponents $t$ and $r$ of variables $u$ and $v$, respectively, in $G'_{i,k}(x,u,v)$.\\

Next, we illustrate how to proceed to compute the configuration index with the example.

\begin{example}

Let us consider again the weighted majority game with coalitions configuration used in Example 3.5. Take again player one, the first coalition, and we will build $G'_{1,1}(x,u,v)$. Following Theorem 3.6:
$$
\begin{array}{rl}
	G'_{1,1}(x,t_2,t_3,t_4,u,v)=&\prod_{{j} \in C_1 \backslash \{1\}}(1+x^{w_{j}}t_{j}u)\prod_{C_l \notin C^1}(1+x^{w(C_l)}[\prod_{j \in C_l}t_{j}]v) \\
&\\
	=& \bigl[(1+x^{w_2}t_2u)(1+x^{w_3}t_3u)\bigr] \bigl[\bigl(1+x^{w(C_2)}t_2t_3v\bigr)\bigl(1+x^{w(C_3)}t_3t_4v\bigr)\bigr]\\
&\\
	=&1+x^{w_2}t_2u+x^{w_3}t_3u+x^{w_2+w_3}t_2t_3u^2+x^{w(C_2)}t_2t_3v+x^{w_2+w(C_2)}t^2_2t_3uv\\
&\\
	+&x^{w_3+w(C_2)}t_2t^2_3uv+x^{w_2+w_3+w(C_2)}t^2_2t^2_3u^2v+x^{w(C_3)}t_3t_4v\\
&\\
	+&x^{w_2+w(C_3)}t_2t_3t_4uv+x^{w_3+w(C_3)}t^2_3t_4uv+x^{w_2+w_3+w(C_3)}t_2t^2_3t_4u^2v\\
&\\
	+&x^{w(C_2)+w(C_3)}t_2t_3^2t_4v^2+x^{w_2+w(C_2)+w(C_3)}t_2^2t_3^2t_4uv^2\\
&\\
	+&x^{w_3+w(C_2)+w(C_3)}t_2t_3^3t_4uv^2+x^{w_2+w_3+w(C_2)+w(C_3)}t_2^2t_3^3t_4u^2v^2.\\
\end{array}
$$
\bigskip

In the last sum, as in Example 3.5, each term represents a coalition consistent with the coalition configuration for $C_1$ that not contains player 1. The terms of the form $t_j$ indicate the players involved and the exponents of $t_j$ indicate the number of times that the corresponding players appear in each subcoalition of the consistent coalition. The exponent of $x$ is the sum of the weights of the subcoalitions. Moreover, the exponent of $u$ is the size of the subset of $C_1$ considered and the exponent of $v$ is the number of coalitions of $C\backslash C_1$.

Now, we multiply each monomial by terms of the form $x^{-(e_{j}-1)w_{j}}$ being $e_{j}$ the exponent of each $t_{j}$ on the coefficient of the monomial. Then, we delete such coeficients and we substitute the weights of each player and the weights of each coalition. In this way, we obtain the generating function of the numbers $\{\nu^{1,1}_{m,r,t}\}_{m \geq 0,\, r\geq 0,\, t\geq 0}$:

\bigskip
\begin{center}
\begin{tabular}{rl}
$G'_{1,1}(x,u,v)=$ & $1+2x^2u+x^4u^2+x^3v+x^3uv+x^4v+2x^4uv$\\
 + & $x^4u^2v+x^5uv+x^5v^2+2x^5uv^2+x^5u^2v^2+x^5u^2v$.\\

\end{tabular}
\end{center}
\bigskip

Finally, we choose the coefficients of the terms whose exponent of $x$ is between $q-w_1=2$ and $q-1=2$, and we have that $\sigma^{1,1}_{0,1}(P,f;C)=2$. Since player 1 is only in one element of the coalition configuration, we have that:

\bigskip
$$
\noindent \begin{array}{rl}
\Phi_{1}(P,f;C)= & \sum\limits_{C_{k} \in C^{1}}\displaystyle{\sum_{r=0}^{c-c^1}}\displaystyle{\sum_{t=0}^{c_k-1}}
\frac{r!(c-r-1)!}{c!}\frac{t!(c_k-t-1)!}{c_k!}\sigma^{1,k}_{r,t}(P,f;C)\\
&\\
= & \displaystyle{\frac{0!(3-0-1)!}{3!}\frac{1!(3-1-1)!}{3!}2=\frac{1}{9}=0.1111}.

\end{array}
$$
\bigskip

Proceeding analogously with the remaining players, we obtain that the configuration index for this game is: $$\Phi(P,f;C)=\bigl(0.1111111,0.2777778,0.4444444,0.1666667\bigr).$$
The same result is obtained by using the code provided in the Appendix. In this case, we introduce the function \texttt{IndexWeightedMajorityCC} and then we write:

\bigskip
\begin{center}
\begin{tabular}{l}

$\texttt{wi<-c(1,2,2,1)}$\\
$\texttt{q<-3}$\\
$\texttt{C<-list(c(1,2,3),c(2,3),c(3,4))}$\\
$\texttt{IndexWeightedMajorityCC(wi,q,C)}\  \Box$\\

\end{tabular}
\end{center}
\bigskip

\end{example}

\section{The case of the European Union of the 28 members (2015)}

In this section, we compute the generalized Banzhaf-Coleman and configuration indices for the different countries which constitute the {\it European Union} (EU). From 1 July 2013, this community comprises 28 members.

The countries that are part of the European Union are indicated in Table 1 and also the number of seats at their disposal in the European Parliament since the Lisbon Treaty became effective.

If we identify each country, in the order given in Table 1, with a player whose weight is the number of seats, and the quota  $q$ = 376 is given by the simple majority, we can model this situation as a weighted majority game.

In addition to this, we introduce two partitions of the players to represent a more realistic situation. For the definition of these partitions, we will follow criteria inspired by the works of Albizuri and Aurrekoetxea (\cite{alb06}) and Courtin (\cite{cou11}). First partition, denoted by $C_1$, will take into account cultural and geographical positions of the states in a way that we will group those located in the vicinity of Mediterranean and Atlantic, which are at the heart of Europe, in the North, and a part of the newest members of the EU.

The other one, $C_2$, will take into account the economic position of the states, and for this we used the list of sovereign states in Europe by GDP (PPP) per capita (lists of countries by gross domestic product at purchasing power parity per capita, the value of all final goods and services produced within a country in 2013, divided by the average population for the same year) provided from the {\it World Bank}.

Thus, we have:
$$C_1=\{ \{1, 2, 6, 10, 11, 18, 20, 24, 25\}, \{3, 4, 8, 13, 16, 17, 21, 23\},\{5, 7, 12, 14, 15, 19, 22, 26, 28\}, \{9, 27\} \}$$
and
$$C_2=\{ \{1, 7, 9, 10, 11, 14, 18, 20, 27, 28\},\{2, 5, 6, 12, 15, 19, 22, 24, 25, 26\},\{3, 4, 8, 13, 16, 17, 21, 23\} \}.$$

As coalition configurations allow players to be in more than one element simultaneously, we will consider our coalition configuration $C=C_1 \cup C_2$.

In Table 1, we show both power indices studied in this paper for this case, with additional indication for each state, of the coalitions to which it belongs (second and third columns) in each of the two partitions defined on the set of players.

\begin{table}[h!]
{\small
\begin{center}
\begin{tabular}{lccccc}
\hline
{\small State}      &{\small $C_1$}&{\small$C_2$}&{\small Seats}&{\small Generalized}&{\small Configuration}  \\
                    &        &     &       &{\small Banzhaf-Coleman}  &  {\small index} \\
                    &        &     &       &   {\small index}         &         \\ \hline
1. Austria          & 1      & 1   &   18  &   0.039794922    & 0.025046296                \\
2. Belgium          & 1      & 2   &   21  &   0.045776367    & 0.030555556                \\
3. Bulgaria         & 2      & 3   &   17  &    0.037109375   & 0.023750000                \\
4. Croatia          & 2      & 3   &   11  &    0.024902344   & 0.016289683               \\
5. Cyprus           & 3      & 2   &   6   &   0.009277344    & 0.009834656                \\
6. Czech Republic   & 1      & 2   &   21  &    0.045776367   & 0.030555556                \\
7. Denmark          & 3      & 1   &   13  &   0.019287109    & 0.009338624                \\
8. Estonia          & 2      & 3   &   6   &   0.012207031    & 0.008710317                \\
9. Finland          & 4      & 1   &   13  &    0.012207031   & 0.007037037                \\
10. France          & 1      & 1   &   74  &   0.157226562    & 0.101243386                \\
11. Germany         & 1      & 1   &   96  &  0.208007812     & 0.138511905                \\
12. Greece          & 3      & 2   &   21  &  0.033325195     & 0.027996032                \\
13. Hungary         & 2      & 3   &   21  &  0.044433594     & 0.029345238                \\
14. Ireland         & 3      & 1   &   11  &  0.014892578     & 0.007645503                \\
15. Italy           & 3      & 2   &   73  &   0.100097656    & 0.094914021                \\
16. Latvia          & 2      & 3   &   8   &  0.014648438     & 0.009900794                \\
17. Lithuania       & 2      & 3   &   11  &  0.024902344     & 0.016289683                \\
18. Luxembourg      & 1      & 1   &   6   & 0.015380859      & 0.013068783                \\
19. Malta           & 3      & 2   &   6   & 0.009277344      & 0.009834656                \\
20. Netherlands     & 1      & 1   &   26  & 0.060058594      & 0.035350529                \\
21. Poland          & 2      & 3   &   51  & 0.113281250      & 0.077757937                \\
22. Portugal        & 3      & 2   &   21  &  0.033325195     & 0.027996032                \\
23. Romania         & 2      & 3   &   32  & 0.070312500      & 0.051289683                \\
24. Slovakia        & 1      & 2   &   13  & 0.024047852      & 0.016117725                \\
25. Slovenia        & 1      & 2   &    8  & 0.021118164      &  0.015185185               \\
26. Spain           & 3      & 2   &   54  & 0.071166992      & 0.065939153                \\
27. Sweden          & 4      & 1   &   20  & 0.020996094      & 0.012513228               \\
28. United Kingdon  & 3      & 1   &   73  &  0.130371094     & 0.087982804                \\\hline
\end{tabular}
\end{center}
\caption{{\protect\small {Power indices for EU 2015 with coalition configuration.}}}
}
\label{Indexs}
\end{table}

Looking at the results of the generalized index of Banzhaf-Coleman, among other findings, we note that the two highest values correspond to France and Germany, which are also the states with the highest number of
deputies. They also belong to the same geographical area and share position economically. The latter also occurs with Cyprus and Malta, which also belong to the states with fewer deputies, and the two states with lower value. Two new members, Croatia and Lithuania, have the same power as they have the same number of deputies and similar economic position. Greece, Portugal, Belgium, and the Czech Republic have the same number of deputies and similar economic position in terms of GDP but the index assigned a lower value to the first two, which is motivated because we have included them in a different geographic area.

If we now analyze the results with the configuration index, we observed similar relationships but not identical, since, as we have seen above, this index takes into account the size of the significant coalitions for each player. According to this index, the country with less power is Finland, which is not the state that least deputies have, but we have included it in a geographic area with only two states. Comparing the results for Italy and the United Kingdom, the former has more power according to the configuration index and less power with the generalized index of Banzhaf-Coleman; both have the same deputies and they are in the same region, but they are in different coalitions depending on the economic level. Estonia and Ireland are also sorted in one way or another according to the index at hand; in this case, the number of deputies, the geographical area, and the economic position are different.

\vspace*{0.5cm}

\section{Conclusions}

Since Albizuri and Aurrekoetxea (\cite{alb06}) and Albizuri, Aurrekoetxea, and Zarzuelo (\cite{alb062}) introduced the generalized Banzhaf-Coleman index and the configuration index, respectively, for simple games with coalition configurations, as a complement to the definition and the study of their properties done in these papers, an open question was the efficient computation of these indices by means of the so-called generating functions, in line with similar studies concerning, among other, to the well-known Shapley-Shubik and Banzhaf-Coleman indices. In this work the task is done, and the construction of computational algorithms, whose suitability is mathematically proved, is complemented with the computer programming using a free software tool such as \texttt{R} and the presentation of examples, both purely numerical in order to explain the algorithms, as taken from real life that show the scope of the model considered and the algorithms introduced. We believe it is worth noting that the model of cooperative games with coalition configurations studied in this work, is being of the interest of others, today, as Albizuri and Vidal-Puga (\cite{alb15}) or Andjiga and Courtin (\cite{and13}), among others. This means that the algorithms introduced in this article may be under consideration or adaptation in future works.

\section*{Appendix. \texttt{R} code}

We have implemented two functions in the open language \texttt{R} that calculate the generalized Banzhaf-Coleman and the configuration indices for weighted majority games with coalition configuration using our algorithms created from its generating functions.

The reason for using \texttt{R} is that it is a free use software tool and it is increasingly widespread in the field of statistics, operations research, and its applications. \texttt{R} allows to create lists from a dataset and it has several commands that allow its manipulation easily and quickly. The functions we created give the result providing as input the vector of weights of the players in the weighted majority game, the quota, and the coalitions that constitute the configuration.\\

The \texttt{R} procedure in order to calculate the generalized Banzhaf-Coleman index for weighted majority games:\\

{\small
\begin{verbatim}
BanzhafWeightedMajorityCC<-function(wi,q,C){
p<-length(wi)
c<-length(C)
Coali<-matrix(0,nrow=p,ncol=c)
SigB<-list()
Banzhaf<-numeric(p)
for(i in 1:p)
{for(g in 1:c){Coali[i,g]<-sum(C[[g]]==i)}
SigB[[i]]<-numeric(sum(Coali[i,]==1))
NOi<-which(Coali[i,]==0)
for(l in which(Coali[i,]==1)){
new_wi<-wi[C[[l]][-which(i==C[[l]])]]
var<-rep(1,length(new_wi))
var_t<-list()
for(play in C[[l]][-which(i==C[[l]])]){
var_t<-c(var_t,play)}
for(v in NOi){
new_wi<-c(new_wi,sum(wi[C[[v]]]))
var<-c(var,rep(2,length(sum(wi[C[[v]]]))))
var_t<-c(var_t,list(C[[v]]))}
exponents<-0
z<-matrix(0,nrow=2)
for(j in 1:length(new_wi)){
friends<-t(combn(var,j))
combinations<-t(combn(length(var),j))
for(k in 1:dim(combinations)[[1]]){
players<-matrix(0,ncol=p,nrow=dim(combinations)[1])
for(play in 1:j){
players[k,var_t[combinations[k,play]][[1]]]<-players[k,var_t[combinations[k,play]][[1]]]+1}
players<-players-1
players[players<0]<-0
exponents<-c(exponents,sum(new_wi[combinations[k,]])-sum(wi*players[k,]))
z<-cbind(z,matrix(c(sum(friends[k,]==1),sum(friends[k,]==2)),nrow=2))}}
z<-z[,which(exponents>=q-wi[i] & exponents<=q-1)]
exponents<-exponents[which(exponents>=q-wi[i] & exponents<=q-1)]
if(is.matrix(z)){
SigB[[i]][l]<-dim(z)[2]
}else{SigB[[i]][l]<-1}
Banzhaf[i]<-Banzhaf[i]+SigB[[i]][l]/(2^(c+length(C[[l]])-2))}}
return(Banzhaf)}
\end{verbatim}
}

The \texttt{R} procedure in order to calculate the configuration index for weighted majority games:\\

\small{
\begin{verbatim}
IndexWeightedMajorityCC<-function(wi,q,C){
p<-length(wi)
c<-length(C)
Coali<-matrix(0,nrow=p,ncol=c)
SigS<-list()
Index<-numeric(p)
for(i in 1:p){
for(g in 1:c){Coali[i,g]<-sum(C[[g]]==i)}
SigS[[i]]<-list()
NOi<-which(Coali[i,]==0)
for(l in which(Coali[i,]==1)){
new_wi<-wi[C[[l]][-which(i==C[[l]])]]
var<-rep(1,length(new_wi))
var_t<-list()
for(play in C[[l]][-which(i==C[[l]])]){
var_t<-c(var_t,play)}
for(v in NOi){
new_wi<-c(new_wi,sum(wi[C[[v]]]))
var<-c(var,rep(2,length(sum(wi[C[[v]]]))))
var_t<-c(var_t,list(C[[v]]))}
exponents<-0
z<-matrix(0,nrow=2)
for(j in 1:length(new_wi)){
friends<-t(combn(var,j))
combinations<-t(combn(length(var),j))
for(k in 1:dim(combinations)[[1]]){
players<-matrix(0,ncol=p,nrow=dim(combinations)[1])
for(play in 1:j){
players[k,var_t[combinations[k,play]][[1]]]<-players[k,var_t[combinations[k,play]][[1]]]+1}
players<-players-1
players[players<0]<-0
exponents<-c(exponents,sum(new_wi[combinations[k,]])-sum(wi*players[k,]))
z<-cbind(z,matrix(c(sum(friends[k,]==1),sum(friends[k,]==2)),nrow=2))}}
z<-z[,which(exponents>=q-wi[i] & exponents<=q-1)]
exponents<-exponents[which(exponents>=q-wi[i] & exponents<=q-1)]
SigS[[i]][[l]]<-matrix(0,nrow=c-sum(Coali[i,]==1)+1,ncol=length(C[[l]]))
row.names(SigS[[i]][[l]])<-0:(c-sum(Coali[i,]==1))
colnames(SigS[[i]][[l]])<-0:(length(C[[l]])-1)
for(r in 0:(c-sum(Coali[i,]==1))){
for(t in 0:(length(C[[l]])-1)){
if(is.matrix(z)){SigS[[i]][[l]][r+1,t+1]<-sum(z[2,]==r & z[1,]==t)
}else{SigS[[i]][[l]][r+1,t+1]<-sum(z[2]==r & z[1]==t)}
Index[i]<-Index[i]+(factorial(r)*factorial(c-r-1)*factorial(t)*factorial(length(C[[l]])-t-1)
*SigS[[i]][[l]][r+1,t+1])/(factorial(c)*factorial(length(C[[l]])))}}}}
return(Index)}
\end{verbatim}
}

\section*{Acknowledgements}

Financial support from {\it Ministerio de Ciencia e Innovaci\'on} and {\it Ministerio de Econom\'{\i}a y Competitividad} of Spain, and FEDER through grants MTM2011-27731-C03 and MTM2014-53395-C3-2-P, is gratefully acknowledged.

\end{document}